\documentclass[a4paper,10pt]{article}

\usepackage[english]{babel}
\usepackage{amssymb}
\usepackage{amsbsy}
\usepackage{amsmath}
\usepackage{epsfig}
\usepackage{graphicx}
\usepackage{psfrag}
\usepackage{pstricks}
\usepackage{a4}

\newcommand{\R}{{\mathbb R}}

\usepackage{amsfonts}

\newcommand{\x}{{\mathsf x}}
\newcommand{\y}{{\mathsf y}}
\newcommand{\z}{{\mathsf z}}

\numberwithin{equation}{section}

\newtheorem{theorem}{Theorem}

\newtheorem{remark}{Remark}

\newtheorem{ex}{Example}

 1   
 1  

\def\cal#1{\mathcal{#1}}

\def\e{\epsilon}

\def\fpd#1#2{\frac{\partial #1}{\partial #2}}
\def\R{{\rm I\kern-.20em R}}
\def\F{{\rm I\kern-.20em F}}

\def\ovl#1{\overline{#1}}

\def\half{\mbox{$\frac{1}{2}$}}

\def\bea{\begin{eqnarray*}}
\def\eea{\end{eqnarray*}}

\def\fpdt#1#2#3{\frac{\partial^2 #1}{\partial #2 \partial #3}}

\def\F{{\bf F}}
\def\ts{\textstyle}

\begin{document}
\title{Dynamics of the Tippe Top via Routhian Reduction}
\author{ M.C.~Ciocci$^1$ and B. Langerock$^2$ \\
\\
\parbox{12cm}{\small
$^1$ Department of Mathematics, Imperial College London, London SW7 2AZ, UK\\
\small
$^2$ Sint-Lucas school
of Architecture, Hogeschool voor Wetenschap \& Kunst, B-9000
Ghent, Belgium}
}
\maketitle

\begin{abstract}
We consider a tippe top modeled as an eccentric sphere, spinning on a horizontal table and subject to a sliding friction.
Ignoring translational effects, we show that the system is reducible using a Routhian reduction technique. The reduced system is a two dimensional system of second order differential equations, that allows an elegant and compact way to retrieve the classification of tippe tops in six groups according to the existence and stability type of the steady states.
\end{abstract}

\section{Introduction} \label{intro}
The tippe top is a spinning top, consisting of a section of a sphere
fitted with a short, cylindrical rod (\emph{the stem}). One typically sets off the top by making it spin with the stem upwards, which we call, from now on the {\em initial spin} of the top. When the top is spun on a table, it will turn the stem down towards the table. When the stem touches the table, the top overturns and starts
spinning on the stem. The overturning motion as
we shall see is a transition from an unstable (relative) equilibrium
to a stable one. Experimentally, it is known that such a transition only occurs when the spin speed exceeds a certain critical value.

Let us make things more precise and first describe in some detail the model used throughout this paper. The tippe top is assumed to be a spherical rigid body with its center of mass $\epsilon$ off the geometrical center of the sphere with radius $R$. In addition, we assume that this eccentric sphere has a mass distribution which is axially symmetric about the axis through the center of mass $O$ and the geometrical center $\cal C$ of the sphere. The tippe top is subject to a holonomic constraint since we only consider motions of the tippe top on a horizontal plane. An important observation is that the system has an integral of motion, regardless of the model of the friction force that is used. This integral is called the {\em Jellet} $J$ and is proportional to the initial spin $n_0$ (equation~\eqref{Jellet}). For the sake of completeness we mention here that a nonholonomic model (rolling without slipping) for the tippe top is not appropriate, since the equations of motion then allow two more integrals of motion (the energy and the Routh-integral) that prohibit the typical turning motion, cf.\ eg.\ \cite{CBJB,Gray}. We will choose the friction force to be linearly dependent on the slip velocity of the contact point $Q$ of the top and the horizontal plane.

In earlier works \cite{scheck,RW} conditions of stability for asymptotic states of the tippe top were retrieved by the Lyapunov function method starting with the Newton equations of motion. We take a different approach to the problem.
Our goal is to describe the behaviour of the tippe top in terms of Lagrangian variables and thereby classify its asymptotic motions as function of non-dimensional (physical) parameters $A/C$, the  \emph{inertia ratio} and $\epsilon/R$, the \emph{eccentricity} of the sphere and of the value the Jellet $J$. As main novelty, we show that by ignoring the translational velocity of the center of mass in the friction law, the Lagrangian formulation of the dynamics of the tippe top allows a restriction to a system on $SO(3)$ which is amenable to Routhian reduction~\cite{pars}. The reduced system allows a rather simple stability analysis of the relative equilibria based on the relationship of the value of the Jellet's integral $J$ and the tumbling angle $\theta$, see below.
Through the Routhian reduction we retrieve a simple stability criterium which leads to the same conclusions as in \cite{scheck,RW}.

The hypothesis of neglecting translational effects is similar to the one made by Bou-Rabee et al.~\cite{Romero}. They motivated it by reasoning that for all possible asymptotic states (the relative equilibria of the system) the velocity of the center of mass is zero \cite{scheck}, and, for this reason, in a neighborhood of these solutions the translational friction can be neglected.
In~\cite{RW} the dynamics of the spherical tippe top with small friction has been studied without such an approximation and a full analysis of the asymptotic long term dynamics of the system is given in terms of the Jellet and the eccentricity of the sphere and inertia ratio. Remarkable is the fact that the stability results obtained by means of the Routhian reduction procedure fully coincides with the results of~\cite{RW}. This justifies a-posteriori the approximation assumption in the friction force and shows that it accurately describes the behavior of the tippe top.

The main result of this paper is summarized in the following theorem, also compare with \cite{CBJB,RW}.

\begin{theorem}\label{mainthm}
In the approximation of negligible translational effects, a spinning eccentric sphere on an horizontal (perfectly hard) surface subject to a sliding friction is reducible with a Routhian reduction procedure~\cite{pars}. The relative equilibria of the reduced system are precisely the steady states of the original system. They are purely rolling solutions and except for the trivial state of rest, they are of three types:
\begin{itemize}
\item[(i)] (non-inverted) vertically spinning top with center of mass straight below the geometric center;
\item[(ii)] (inverted) vertically spinning top with center of mass straight above the geometric center;
\item[(iii)] intermediate spinning top, the top precesses about a vertical while spinning about its axle and rolling over the plane without gliding.
\end{itemize}
The existence and stability type of these relative equilibria only depend on the inertia ratio $\frac{A}{C}$, the eccentricity of the sphere $\frac{\epsilon}{R}$ and the Jellet invariant $J$. In particular, six regimes are identified in terms of the Jellet invariant and only three exhibit the `tipping' behavior.
\end{theorem}
It turns out that the vertical states always exist, and intermediate states may branch off from them. Qualitative bifurcation diagrams corresponding to the possible different regimes are sketched in Fig.~\ref{fig:bifdia}.

\newcommand{\schaal}{.5}
\begin{figure}[ht]
\centering \includegraphics[scale=\schaal]{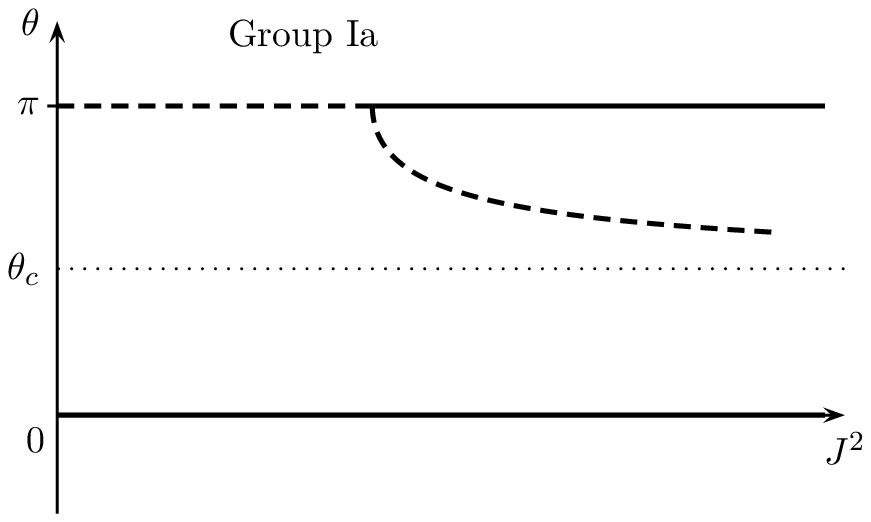}\hspace{.5cm}
\includegraphics[scale=\schaal]{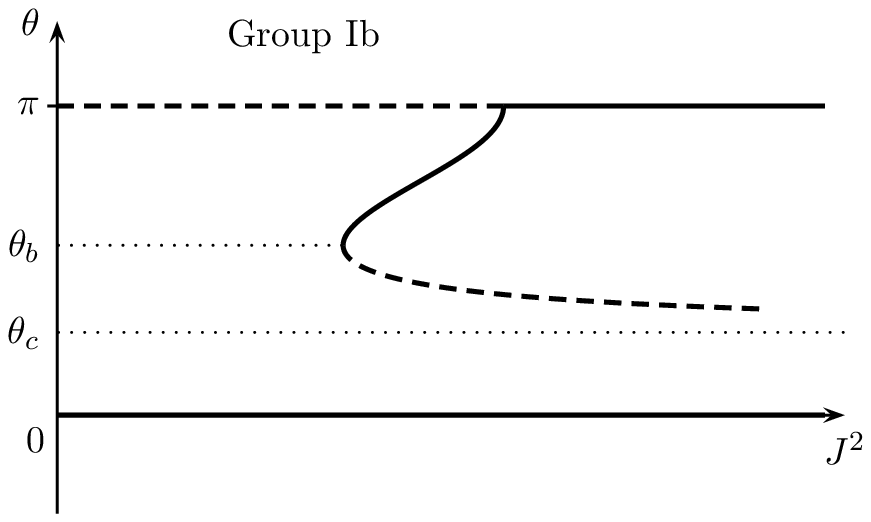}\\
\centering \includegraphics[scale=\schaal]{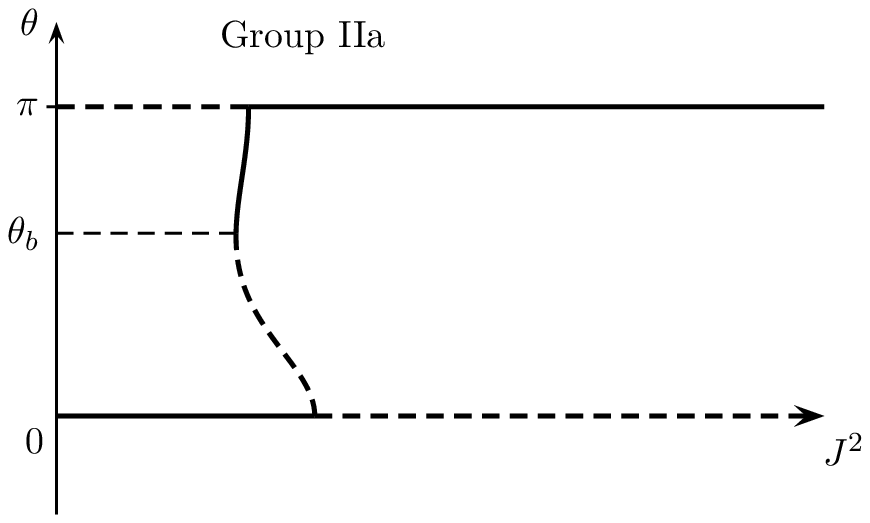}\hspace{.5cm}
\includegraphics[scale=\schaal]{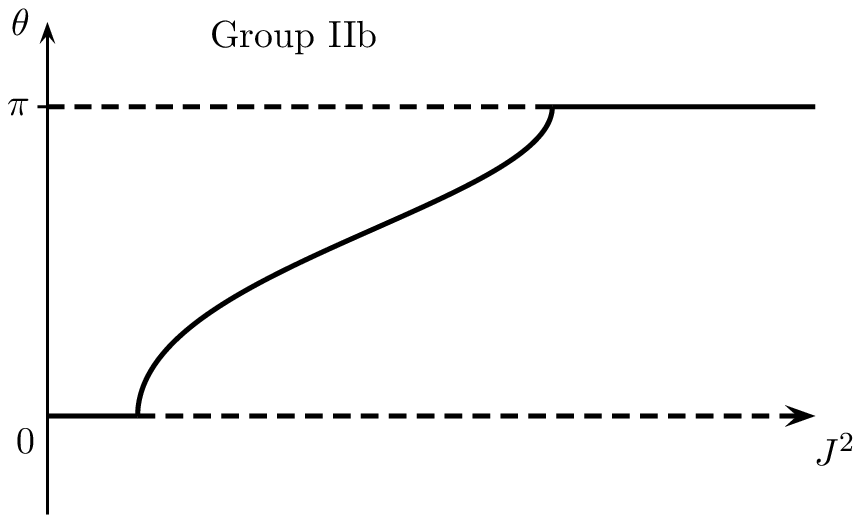}\\
\centering \includegraphics[scale=\schaal]{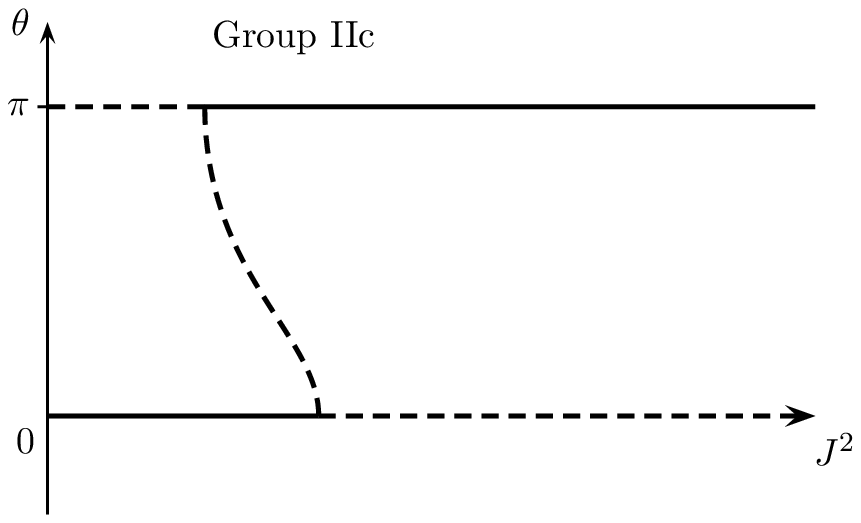}\hspace{.5cm}
\includegraphics[scale=\schaal]{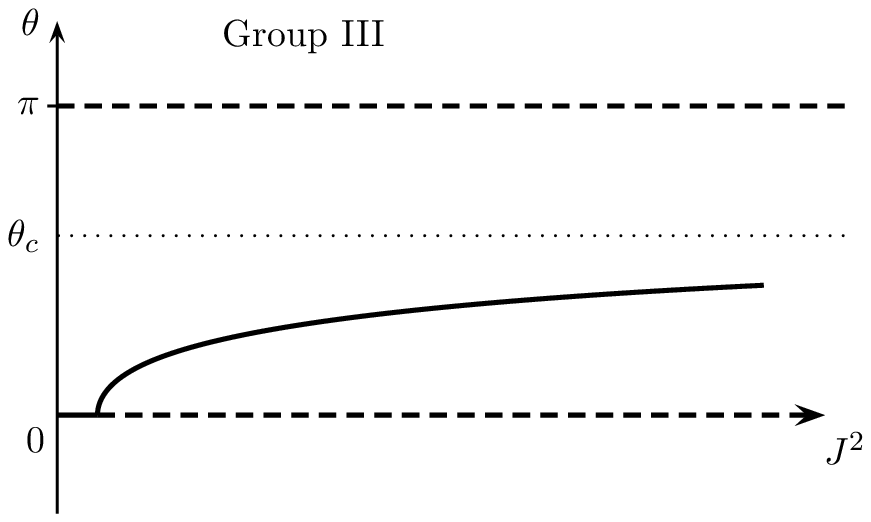}
\caption{\footnotesize Bifurcation diagrams of relative equilibria in function of the Jellet invariant $J$. Solid black branches correspond to stable relative equilibria, while dashed black branches correspond to unstable ones. The vertical states
$\theta=0$ and $\theta=\pi$ always exist, from which intermediate states (with $0<\theta<\pi$) may branch off.}
\label{fig:bifdia}
\end{figure}
The reader is referred to Sec.~\ref{Sec:Stability} -- Sec.~\ref{Sec:stability} for the
details and specification of the parameter ranges. There are three main groups Group I, II and III determined by existence properties of intermediate equilibrium states (similar to the subdivision from~\cite{watanabe}). Tippe tops of Group I may admit intermediate states for $\theta>\theta_c$, where $\theta$ denotes the angle between the vertical and the symmetry axis of the top. Two subgroups are distinguished  according to change in stability type of the intermediate states.
Relevant is that for tippe tops belonging to this class the
non-inverted position is always stable, so they never flip  however
large the initial spin when launched under an angle $\theta$ close
to $0$.
Tippe tops of Group II may admit intermediate states for all
$\theta$, and they show complete inversion when the initial spin is large enough. Tops of Group III tend to flip over up to a certain angle
$\theta_c<\pi$ when spun rapidly enough. Since stability results are often
in the literature expressed using the `initial' spin $n_0$ of the
tippe top, one can read the $J^2$ in the figures as $n_0^2$. We anticipate to further results by noting that the instability inequalities are independent of the friction coefficient unless it is zero.

\medskip
\noindent
The structure of the paper is as follows. In
Section~\ref{sec:equat} our model is described and the equation of motion are given according to the Lagrangian formalism. After introducing the Jellet integral of motion, the Routhian reduction is performed in Section~\ref{sec:Routhreduction}. The steady states of the system are then calculated and their stability type is determined, yielding a tippe top classification in six groups which is summarized in Section~\ref{classifi}.

\section{Equations of motion}\label{sec:equat}
\begin{figure}[htb]
\begin{center}
\scalebox{.5}{ \psfrag{t_1}{$z$} \psfrag{t_2}{$x$} \psfrag{t_3}{}
\psfrag{t_4}{$\z$} \psfrag{t_5}{$\psi$} \psfrag{t_6}{$\varphi$}
\psfrag{t_7}{$\theta$} \psfrag{t_8}{$\mathcal C$} \psfrag{t_9}{$O$}
\psfrag{e}{$\epsilon$} \psfrag{R}{$R$} \psfrag{t_10}{$Q$}
\psfrag{t_11}{$h(\theta)=R-\epsilon \cos(\theta)$} \psfrag{t_12}{$\bf G$} \psfrag{t_13}{$\bf{R}_f$}
\psfrag{t_14}{$\Pi$}
\includegraphics[width=11cm]{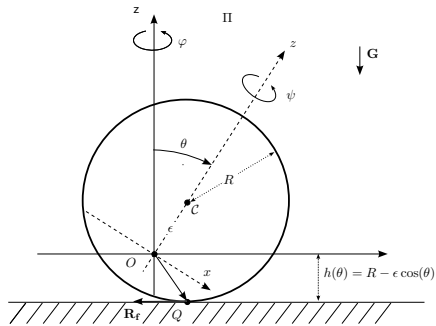}}
\caption{\footnotesize Eccentric sphere version of the tippe top. $R$ is the
radius of the sphere, the center of mass $O$ is off center by
$\epsilon$. The top spins on a horizontal table with point of
contact $Q$. The axis of symmetry is $Oz$ and the vertical axis is
$O\z$, they define a plane $\Pi$ (containing $\vec{OQ}$) which
precesses about $O\z$  with angular velocity $\dot \varphi$. The
height of $O$ above the table is $h(\theta)$.
\label{tiptopgeometry}}
\end{center}
\end{figure}

As is mentioned in the introduction, we consider the \emph{eccentric sphere} model of such a top, see Fig.~\ref{tiptopgeometry}. That is, we consider a sphere with radius $R$ whose mass
distribution is axially symmetric but not spherically symmetric, so
that the center of mass and the geometric center do not coincide.
The line joining the center of mass and the geometrical center is an
axis of inertial symmetry, that is, in the plane perpendicular to
this axis the inertia tensor of the sphere has two equal
principal moments of inertia $A=B$. The inertia moment along the axis of symmetry is denoted by $C$ and the total mass of the sphere is $m$.

The eccentricity $\epsilon$ is the distance between the center of mass $O$ and the geometric center
$\cal C$ of the sphere, with $0<\epsilon<R$. The point $Q$ is the
point of contact with the plane of support.
We assume that an inertial (laboratory) frame $M\x\y\z$ is chosen, where $M$ is some point on
the table and the $\z$-axis is the vertical. Let us denote the unit vectors along the axis of the reference frame
$O\x,O\y,O\z$ fixed to the body by respectively $e_\x,e_\y,e_\z$. The coordinates of the center of mass are denoted by ${\bf r}_O=(\x,\y,\z)_{M\x\y\z}$. A second reference frame is denoted by $Oxyz$, and is defined in such a way that its third axis is precisely the symmetry axis of the top and the $y$-axis is perpendicular to the plane $\Pi$ through the $\z$- and $z$-axes (see Fig.~\ref{tiptopgeometry}). Again we denote the unit vectors along the axis of this reference frame by $e_x,e_y,e_z$\footnote{In~\cite{Romero} the origin of the reference system attached to the body is in
the center of the sphere $\cal C$, and not in the center of mass $O$.}.

Let $(\theta,\varphi,\psi)$ be the Euler angles of the
body with respect to the inertial frame, Fig.~\ref{tiptopgeometry}, chosen in such a way that (i) the vertical plane $\Pi$ is
inclined at $\varphi$ to the fixed vertical plane $\x\z$ and
precesses with angular velocity $\dot \varphi$ around the vertical
$O\z$; (ii) the angle $\theta$ is the angle between the vertical
$O\z$ and the axle $O\z$ of the top; $\dot{\theta}$ causes the nodding
(nutation) of the axle in the vertical plane $\Pi$; and (iii) the angle $\psi$
orients the body with respect to the fixed-body frame, $\dot
\psi$ is the spin about the axle.

As it was pointed out before, the tippe top is constrained to move
on a horizontal plane. This holonomic constraint is expressed by
$\z=R-\e\cos\theta=h(\theta)$. We assume throughout the paper that the only forces acting on the sphere are gravity ${\bf G}=-mg e_\z$ and a friction force ${\bf{F}}$ exerted at the point of contact $Q$ of the sphere with the plane. It is now immediate to write down the Lagrangian for the tippe top:
\begin{multline}
L=\half \left(m(\dot \x^2 + \dot \y^2)+ (\e^2m\sin^2\theta+A)\dot
\theta^2  + A\sin^2\theta\dot\varphi^2+ C(\dot\psi
+\dot\varphi\cos\theta)^2\right) \\   - mg(R-\e\cos\theta),
\end{multline}
where $g$ is earth acceleration.

This function is defined on the tangent space of the configuration manifold $M=\R^2\times
SO(3)$. In order to obtain the equations of motion, it only remains to define a suitable friction force ${\bf{F}}$. The Lagrangian equations of motion for the tippe top then read:
\begin{equation}\label{thesystem}
\frac{d}{dt}\left(\frac{\partial L}{\partial \dot q^i}
 \right)-\frac{\partial L}{\partial q^i}=Q^{F}_i,
\end{equation}where $q^i$ represents one of the coordinates $(\x,\y,\varphi,\theta,\psi)$, and $Q^{F}=Q^{F}_idq^i$ is a one-form on $M$ representing the generalized force moment of the friction force at the point of contact. It is defined by, with ${\bf F}={\bf R}_f+R_n e_\z$ the orthogonal decomposition of ${\bf F}$:
\[
Q^{F}={\bf{R}}_f\cdot e_\x d\x + {\bf{R}}_f\cdot e_\y d\y+ ({\bf{q}}
\times {\bf{R}}_f)\cdot (e_y d\theta + e_\z d\varphi +
e_zd\psi).
\]

\paragraph{Modeling the friction force}
One typically models the friction force ${\bf F}={\bf R}_f+R_ne_\z$ to be proportional to the slip velocity of the point of contact ${\bf v}_Q$. We denote by
${{\bf{R}}_n}=R_n {e}_\z$ the normal reaction of the floor at
$Q$, which is of order $mg$, and ${\bf R}_f= F_X {e}_\x+F_Y
{e}_\y$ is the (sliding) friction which opposes the slipping
motion of the body. The fact that the sliding friction opposes the slipping motion is
expressed by ${\bf R}_{f}\cdot {\bf{v}}_Q \leq 0$.
We adopt a viscous friction law~\cite{CBJB,Levi,MofShi,watanabe} and assume that
\begin{equation}\label{viscfric}
{\mathbf R}_f=-\mu R_n {\bf{v}}_Q.
\end{equation}
Here $\mu$ is a coefficient of friction with the dimension of
(velocity)$^{-1}$. It now takes a few tedious computations to arrive to the coordinate expressions for the force moments of the friction force.
The coordinates of the point of contact $\vec{OQ}:={\bf q}$ are $Q =(x_Q,0,z_Q)_{Oxyz}=(R\sin\theta,0,\e-R\cos\theta)_{Oxyz}.$
The velocity of the point of contact $Q$ equals
\begin{equation}
{\mathbf{v}}_Q={{\bf{v}}_O} + {\boldsymbol{\omega}}\times {\mathbf{q}},
\end{equation}
${{\mathbf{v}}_O}=(\dot \x,\dot \y,h^\prime(\theta)\dot\theta)_{M\x\y\z}$ is the velocity of the
center of mass. A coordinate expression for the angular velocity is given by
\begin{equation} \label{omega}
{\boldsymbol{\omega}}=  -\dot{\varphi} \sin(\theta) e_x \; + \; \dot{\theta} e_y \; + \; {n} e_z, \quad \mbox{where} \quad {n}:=\dot{\psi}+\dot{\varphi}\cos(\theta)
\end{equation}
$n$ is the spin (that is, the component of $\boldsymbol{\omega}$ about $Oz$). The generalized force moments now read:
\begin{eqnarray*}
Q_\x&=& -\mu R_n (\dot \x-\sin\varphi\ \dot
\theta(R-\e\cos\theta)+\cos\varphi\sin\theta(R\dot\psi+\e\dot
\varphi))\\ Q_\y&=&-\mu R_n(\dot \y+\cos\varphi\ \dot
\theta(R-\e\cos\theta)+\sin\varphi\sin\theta(R\dot\psi+\e\dot \varphi))\\
Q_\theta&=& -\mu R_n(R-\e\cos\theta)(\cos\varphi\dot \y -\sin\varphi
\dot \x
+(R-\e\cos\theta)\dot\theta)\\
Q_\varphi&=& -\mu R_n\e \sin\theta(\cos\varphi\dot \x+\sin\varphi\dot \y+\sin\theta (\e\dot \varphi+R\dot\psi))\\
Q_\psi&=& -\mu R_n R\sin\theta(\cos\varphi\dot \x+ \sin\varphi \dot
\y +\sin\theta(R\dot \psi +\e\dot\varphi)).
\end{eqnarray*}
For the sake of completeness we write an explicit expression for the normal component of the reaction force $R_n$, which can be determined from Newton's law for the center of mass of the sphere:
\begin{equation}\label{reactie}
R_n(\theta,\varphi,\dot\theta,\dot\varphi,\dot\psi,\dot\x,\dot\y)=\frac{g+\dot
\theta^2 h^{\prime
\prime}+h^{\prime}\dot\varphi\sin(\theta)(\dot\varphi\cos(\theta)-C
({\dot\psi
+\dot\varphi\cos\theta})/A)}{1/m+h^{\prime}/A[-h\mu(\sin\varphi\dot\x-\cos\varphi\dot\y-\dot\theta
h)+ h^{\prime}]}
\end{equation}

\noindent To conclude this section, we briefly discuss other models for the tippe top. The eccentric sphere model does not accurately model the contact effects of the tippe top stem. However,
it does describe the fundamental phenomenon of the over-turning. For the sake of
completeness we remark that as soon as the top rises to spin on its
stem, it behaves as a `normal' spinning top with rounded peg, we
refer to~\cite[Chapter 6]{book} for a satisfactory introduction
to this topic. Note that our model also does not describe the peculiar `Hycaro'
tippe tops by Prof.\ T.\ Tokieda~\cite{T^2}, which need a non-axisymmetric asymmetric mass distribution. These tippe tops have a
`preferred direction' meaning that the top would flip over only if
spun in the preferred direction with a certain initial spin, and, no
matter what the initial spin is, it would just continue rotating
around the rest position when spun the other way round. We refer to~\cite{Duitsers} for a detailed analysis when elasticity
properties of the horizontal surface and tippe top are taken into account. Their model allows for jumps of the tippe top on the horizontal surface (that we assumed to be rigid). Finally, we mention that we do not debate over the issue of whether transitions sliding-rolling and rolling-sliding occur in the motion of an eccentric sphere on a flat surface. We chose to concentrate on the sliding model only, because we were interested in capturing the `overturning' phenomenon which cannot occur under the non-holonomic constraint of pure rolling. We refer the interested reader to \cite{CBJB,Gray,Levi} for a discussion of the topic and to \cite{cush,mamaev,Tor} for an analysis of the motion of the rolling eccentric sphere also called the Routh's sphere.

\subsection{Constants of motion: the Jellet invariant }
It was first shown by Jellet~\cite{Jellet} by an approximate argument, and later proved by Routh~\cite{routh} that the system, even if dissipative, has a conserved quantity:
\begin{equation}\label{Jellet}
J=-{\bf{L}}\cdot {\bf{q}}= {\rm const},
\end{equation}
where ${\bf L}$ is the angular momentum of the tippe top about the center of mass. We prove this by using Euler equations which govern the evolution of the angular momentum $\dot{{\bf {L}}}={\bf{q}}\times {\bf{F}}.$
The total time derivative of $J$ then becomes:
\begin{align*}
  \dot J&=-\dot{\mathbf L}\cdot {\bf q} -{\mathbf L}\cdot \dot {\mathbf q}\\
  &=0-\left(A{\boldsymbol{\omega}}-(A-C)({\boldsymbol{\omega}}\cdot e_z)e_z\right)\cdot (\e {\boldsymbol{\omega}}\times e_z)=0\ .
\end{align*}
Straightforward calculations show that Jellet's constant can be written as
\begin{equation}
 J=C{n}(R \cos(\theta)-\epsilon)+A \dot{\varphi} R
\sin^2(\theta)\ .\label{jello}
\end{equation}
We emphasize once more that the Jellet's constant is an
exact constant of motion for the tippe top whether or not there is
slipping and independent of the expression for $\bf F$. As we will explain later, it is this constant that to some
extent controls the motion of the spinning top.  Indeed, it allows a
Routhian reduction procedure (see Sec.~\ref{sec:Routhreduction}),
resulting in relatively simple reduced equations from which we are
able to recover in full detail the stability properties of the
steady states. In the specific case that $\theta=0$, the Jellet is proportional to $n_0=n|_{\theta=0}$, the spin about the $z$-axis. Since one typically sets off the tippe top at an angle $\theta\approx 0$, one can say that the Jellet is proportional to the initial spin $n_0$. Note that the spin at $\theta=\pi$ has an opposite sign to the initial spin $n_0$, meaning that, relative to a body fixed frame, the spin is reversed when the tippe top fulfills a complete inversion.

\medskip
\noindent There is a rotational symmetry for
which the Jellet is the associated first integral. The action of
$S^1$ on $\R^2\times SO(3)$ can be defined as a simultaneous
rotation about $\hat e_3$ over the angle $R\xi$ and about $\hat k$
over the angle $-\e \xi$, where $\xi \in S^1$ (see also~\cite{Romero}). Noether's theorem is applicable in this situation
since the work of the friction force at the point of contact
vanishes under this action.

\medskip
\noindent Note that the total energy of the spinning top is $E=T+V$ is in
general not conserved. Here $T$ is the kinetic energy
with its rotational and translational part, $V=m gh(\theta)$ is the
potential energy. The orbital derivative of $E$ is
\begin{equation}\label{derE}
\frac{d}{dt}E={\bf v}_Q\cdot {\bf R}_{f} \leq 0,
\end{equation}
which is negative semi-definite and vanishes if and only if
${\bf v}_Q$ vanishes. Observe that $E(t)$ decreases monotonically
and hence is a suitable Lyapunov function\footnote{$E(t)$ is
analytical, therefore it is either strictly monotone or a constant.
The energy $E$ is constant only if ${\bf v}_Q=0$. Note that $E$
being Lyapunov~\cite{scheck} implies that the limiting solutions for $t\rightarrow
\infty$  are solutions of constant energy.}, see~\cite{scheck}. From (\ref{derE}) it
follows that dissipation is due to friction.

\subsection{Routhian reduction}\label{sec:Routhreduction}
\newcommand{\etal}{et al.}
It turns out that, if we consider an approximation of the
friction law, the resulting generalized force moments assume a form
that allows us to apply a Routhian reduction procedure, see~\cite{pars} and appendix~\ref{app:routh}. In turn,
using the reduced equations we are able to study in full detail the
stability properties of the tippe top which confirm the results
obtained in~\cite{RW}, and also
recover those of~\cite{Romero,watanabe}.

We now ignore
translational effects in the friction force, i.e. we assume
that all terms in $Q_{\theta}, Q_{\varphi}, Q_{\psi}$ containing
$\dot \x$ and $\dot \y$ are neglected. Typically this approximation is
justified by noting that for all steady states the velocity of the
center of mass is zero, and that in a neighborhood of the steady
states it can be neglected. In our situation, it allows to restrict
ourselves to a system on $SO(3)$ which is reducible using Routh's
procedure. It is easily seen that within this
approximation, if we study the Lagrangian system on $SO(3)$
determined by
\begin{eqnarray*}
L'&=&\half \left( (\e^2m\sin^2\theta+A)\dot \theta^2 +
A\sin^2\theta\dot\varphi^2+ C(\dot\psi
+\dot\varphi\cos\theta)^2\right) \\ && \qquad \quad
-mg(R-\e\cos\theta),\\
Q'&=& -\mu R'_n(R-\e\cos\theta)^2\dot\theta d\theta
 -\mu R'_n\e \sin^2\theta(\e\dot \varphi+R\dot\psi)d\varphi\\ &&\qquad \quad
-\mu R'_n R\sin^2\theta(R\dot \psi +\e\dot\varphi)d\psi,
\end{eqnarray*}
with $R'_n(\theta,\dot \theta,\dot\varphi,\dot\psi) =
R_n(\theta,\varphi,\dot \theta,\dot\varphi,\dot\psi,\dot\x=0,\dot\y=0)$, then we essentially study the entire approximated system. Indeed, any
solution $(\varphi(t),\theta(t),\psi(t))$ to this Lagrangian system
will determine the remaining unknowns $(\x(t),\y(t))$ as solutions
to the following system of time-dependent second order differential
equations:
\begin{equation}
m\ddot \x = Q_\x(\dot\x,\theta(t),\varphi(t),\psi(t)), \quad m\ddot
\y = Q_\y(\dot\y,\theta(t),\varphi(t),\psi(t)).
\end{equation}
Our next step in the reduction procedure is to consider the Lagrangian system $L'$ on $SO(3)$ with generalized
force form $Q'$ and perform a simple coordinate transformation,
determined by
\[
(\theta,\varphi,\psi)\mapsto(\theta,\ovl\varphi=\e\varphi
+R\psi,c=R\varphi-\e\psi).
\]
The Lagrangian $L'$ and the force $Q'$ then become
\begin{align*}
L'&=\ts\half \bigg( (\e^2m\sin^2\theta+A)\dot \theta^2 +
\frac{A\sin^2\theta}{(\e^2+R^2)^2}(\e\dot{\ovl \varphi}+R\dot c)^2+\\
& \qquad \ts \frac{C}{(\e^2+R^2)^2}\big((R+\e\cos\theta)\dot{\ovl\varphi} +
(R\cos\theta-\e)\dot c\big)^2\bigg)- mg(R-\e\cos\theta)\\
 Q'&=
\ts -\mu R'_n\big((R-\e\cos\theta)^2\dot\theta d\theta
+\sin^2\theta\dot{\ovl\varphi} d\ovl\varphi\big),\\
\end{align*}
where it is understood that $R'_n$ is a function of
$\theta,\dot\theta,\dot{\ovl\varphi}$ and $\dot c$. The main reason
for writing $L'$ and $Q'$ in this form is the fact that $c$ is a
cyclic coordinate and that $Q'_c=0$. Indeed, recall that Jellet
integral was associated to the symmetry determined by the vector
field $R\partial_\varphi -\e\partial_\psi$, or in the above
introduced coordinate system by the vector field $\partial_c$. In
particular, we have made the symmetry generator into a coordinate
vector field, and this ensures that we can apply the Routhian
reduction procedure~(\cite{pars} and appendix~\ref{app:routh}), provided the coefficients of $Q'$
do not depend on $c$. This is the case since we neglected the terms
in the velocity of the center of mass. The conserved quantity
associated with the cyclic coordinate is precisely the Jellet
integral:
\begin{align*}
\fpd{L'}{\dot c} &=\ts \frac{R A\sin^2\theta}{(\e^2+R^2)^2}(\e\dot{\ovl
\varphi}+R\dot c)+
\frac{C(R\cos\theta-\e)}{(\e^2+R^2)^2}\big((R+\e\cos\theta)\dot{\ovl\varphi}
+ (R\cos\theta-\e)\dot c\big) \\ & \ts = J/(\e^2+R^2).
\end{align*}
The latter equality implies that
\begin{eqnarray}
\dot c &=&  \frac{J(\e^2+R^2)-(RA\e\sin^2\theta
+C(R\cos\theta-\e)(R+\e\cos\theta))\dot{\ovl\varphi}}{R^2
A\sin^2\theta + C(R\cos\theta-\e)^2}\label{routh1}
\end{eqnarray}
The Routhian reduction procedure defines a Lagrangian system
$${\cal R} = L'-J\dot c/(\e^2+R^2)$$ with two degrees of freedom
$(\theta,\ovl\varphi)$ (here every instance of $\dot c$ in ${\cal
R}$ is replaced using~\eqref{routh1}).  The reduced equations of
motion are then given by
\begin{eqnarray}
\label{reducedeqs} \frac{d}{dt}\left(\fpd{{\cal R}}{\dot
\theta}\right)-\fpd{{\cal R}}{\theta} &=& Q'_\theta =-\mu
R'_n(R-\e\cos\theta)^2\dot\theta \\
\label{reducedeqs2} \frac{d}{dt}\left(\fpd{{\cal
R}}{\dot{\ovl\varphi}}\right)-\fpd{{\cal R}}{\ovl\varphi} &=&
Q'_{\ovl\varphi} =-\mu R'_n\sin^2\theta\dot{\ovl\varphi},
\end{eqnarray}
where it is understood that~\eqref{routh1} is used to eliminate
$\dot c$ in $R'_n$. It takes rather tedious but straightforward
computations to show that $\cal R$ can be written as
$${\cal R} = T_2 + T_1 - W,$$
where
\begin{eqnarray*}
T_2 &=&\ts T_{\theta\theta}(\theta) \dot\theta^2 +
T_{\ovl\varphi\ovl\varphi}(\theta) \dot{\ovl\varphi}^2\\ &=&
\frac12(\e^2m\sin^2\theta +A)\dot\theta^2 + \frac12\frac{AC
\sin^2\theta }{(R^2 A\sin^2\theta +
C(R\cos\theta-\e)^2)}\dot{\ovl\varphi}^2, \\
T_1&=&\ts \frac{J}{(R^2+\e^2)}\frac{( R\e
A\sin^2\theta+C(R+\e\cos\theta)(R\cos\theta-\e))}{R^2
A\sin^2\theta + C(R\cos\theta-\e)^2}\dot{\ovl\varphi},\\
W&=&\ts \frac12\frac{J^2}{R^2 A\sin^2\theta +
C(R\cos\theta-\e)^2}-mg\e\cos\theta .
\end{eqnarray*}

The function $W$ is also called the effective potential.

\begin{remark}
The above defined Routhian function $\cal R$ is not globally defined
on the sphere. In order to provide a globally defined system of
differential equations for the reduced system, we need to extract
the term $T_1$ from the Routhian and consider it as a gyroscopic
force (see e.g.\ \cite{marsdenrouth}).
\end{remark}
\begin{remark}
Observe that the effective potential $W$, obtained through reduction, coincides with the effective energy on a Jellet's level surface as it has been used in \cite{RW}, Sec.~3.
\end{remark}

\medskip
\noindent Note that $\mathcal R$ and $Q'$ are independent of
$\ovl\varphi$. This residual symmetry does not lead to a conserved
quantity (the friction does not vanish for
$\partial_{\ovl\varphi}$). This symmetry is due to the approximation
we carried out in the previous section; it is not present in the
original system $(L,Q)$ or~\eqref{thesystem}. It leads however to a
zero eigenvalue of the linearized system at equilibrium points, see
Section~\ref{Sec:stability}.

\section{Steady states }\label{Sec:Stability}
The equilibria of the reduced Routhian system \eqref{reducedeqs} and \eqref{reducedeqs2} are determined by
\begin{equation}
\dot \theta=0, \dot{\ovl\varphi}=0 \quad \rm{and} \quad
\partial W/\partial \theta =0, \label{crit_amendedpot}
\end{equation}
i.e. they are the critical points of the effective potential (note that if $\dot \theta=\dot{\ovl\varphi}=0$ the components of the force vanish).
Equation (\ref{crit_amendedpot}) is satisfied if (i) $\sin\theta =0 $ or, if (ii)
\begin{multline}\label{define_intstates}
f(J^2,\cos\theta) =\ts
\frac{J^2}{mgCR^2\e}\left(\left(\frac{A}{C}-1\right)\cos\theta+\frac{\e}{R}\right)-
\left(\frac{A}{C}\sin^2\theta +
\left(\cos\theta-\frac{\e}{R}\right)^2\right)^2= 0.
\end{multline}
Solutions to (i) are $\theta=0,\pi$ and give the vertical spinning states. Solutions to (ii) only occur if
\begin{equation}\label{Exist}
\left( \frac{A}{C}-1\right)\cos(\theta)+\frac{\epsilon}{R}>0.
\end{equation}
If this condition is satisfied, solutions to (ii) are the so-called intermediate states. The existence condition only depends on $A/C$ and $\e/R$, and will determine in our classification the three main groups I, II and III, as it was proposed in~\cite{watanabe}. The values for $\theta$ for which (ii) is satisfied will depend on the Jellet, $A/C$ and $\e/R$.

The vertical spinning states correspond to the periodic motion of the tippe top spinning about its axle (which is in vertical position) either in the non-inverted position ($\theta=0$) or inverted position ($\theta=\pi$).
The intermediate states correspond to those relative equilibria in which the tippe top shows in general quasi-periodic motion precessing about a vertical while spinning about its inclined axle rolling over the plane without gliding 
(observe that the intermediate states correspond to the tumbling solution of~\cite{scheck}).

\medskip
The condition (\ref{Exist}) for existence of intermediate states leads to a first classification of tippe tops into three groups.\\
- Group I $[(A/C-1)<-\e/R]$: the tippe tops belonging to this group do not admit
intermediate states in an interval of the form $[0,\theta_c[$, where $\theta_c$ is determined by $(A/C-1)\cos\theta_c
+\e/R =0$.\\
- Group II $[-\e/R<(A/C-1)<\e/R]$: intermediate states may exist for all $\theta\in ]0,\pi[$.\\
- Group III $[(A/C-1)>\e/R]$: tippe tops belonging to this group do not admit intermediate
state in an interval of the form $]\theta_c,\pi[$, where
$\cos\theta_c= (\e/R)/(1-A/C)$.

\medskip \noindent
In the following section we refine this first classification taking into account the stability type of the steady states and their bifurcations, explaining and giving the details of the $J^2$ versus $\theta$ diagrams in Fig.~\ref{fig:bifdia} from the introduction.

We anticipate that the subdivision in subgroups according to a change in stability type of the intermediate states is based on the simple observation that they lie on a curve $f(J^2,\cos\theta) =0$ in the $(J^2,\theta)$-plane, and, denoting by $\partial_2f$ the partial derivative to the second
argument of $f$, a bifurcation point for intermediate states is characterized as the point where
\begin{eqnarray*}
  f(J^2,\cos\theta)&=&0,\\
  \partial_2 f(J^2,\cos\theta)&=& 0.
\end{eqnarray*}
Note that the relation $\partial_2 f(J^2,\cos\theta)= 0$ is essentially the same basic relation (4.27) in \cite{RW}, however its derivation is different. As a consequence we expect the stability results for intermediate states to confirm earlier known facts.

\section{Stability analysis via the reduced equations}\label{Sec:stability}
Determining the (linear) stability of the steady states as given above is an extremely simple task in the reduced setting. Indeed,
let $(\theta_0,\ovl\varphi_0)$ be an equilibrium. The linearized
equations of motion at this relative equilibrium read as
\begin{eqnarray*}
T_{\theta\theta}(\theta_0)\ddot\theta &=& \fpdt{T_1}{\dot{\ovl\varphi}}{\theta}(\theta_0)\dot{\ovl\varphi}
-\fpd{^2W}{\theta^2}(\theta_0)(\theta-\theta_0)
-\mu mg(R-\e\cos\theta_0)^2\dot\theta\\
T_{\ovl\varphi\ovl\varphi}(\theta_0)\ddot{\ovl\varphi}  &=&
-\fpdt{T_1}{\dot{\ovl\varphi}}{\theta}(\theta_0)\dot\theta -\mu mg
\sin^2\theta_0\dot{\ovl\varphi},
\end{eqnarray*}
where we used that $R'_n$ equals $mg$~\eqref{reactie} at the relative
equilibria. It is not hard to show that the
characteristic polynomial of this system is
\begin{eqnarray*} &&p(\lambda)= \lambda\Bigg[\lambda^3 +\mu mg\left(\frac{(R-\e\cos\theta_0)^2}{T_{\theta\theta}(\theta_0)}+ \frac{\sin^2\theta_0}{T_{\ovl\varphi\ovl\varphi}(\theta_0)}\right)\lambda^2 \\  && + \left(\frac{\big(\mu mg(R-\e\cos\theta_0)\sin\theta_0\big)^2 + \left(\fpdt{T_1}{\dot{\ovl\varphi}}{\theta}(\theta_0)\right)^2}{T_{\theta\theta}(\theta_0)T_{\ovl\varphi\ovl\varphi}(\theta_0)}  + \frac{\fpd{^2W}{\theta^2}(\theta_0)}{T_{\theta\theta}(\theta_0)}\right)\lambda
\\ && +\frac{\big(\mu mg \sin^2\theta_0\big) \fpd{^2W}{\theta^2}(\theta_0)}{T_{\theta\theta}(\theta_0)T_{\ovl\varphi\ovl\varphi}(\theta_0)}  \Bigg].\end{eqnarray*} Due to the translational
symmetry in $\ovl\varphi$, one eigenvalue is zero. In Appendix~\ref{app:roots} we show that all remaining eigenvalues have a negative real part, if and
only if $\partial^2W/\partial \theta^2 (\theta_0)
>0$, or if
\begin{multline}\label{stabilityba}
\ts mg\e\cos\theta_0 > \ts \frac{J^2 R^2 C}{(R^2 A\sin^2\theta_0
+ C(R\cos\theta_0-\e)^2)^2}\ts\bigg(-(A/C-1)\sin^2\theta_0\\
\ts + \frac{\mathcal{B}}{C}\cos\theta_0
-4\frac{\mathcal{B}^2}{C^2}\frac{\sin^2\theta_0}{(A/C)\sin^2\theta_0
+ (\cos\theta_0-\e/R)^2}\bigg),
\end{multline}
with $\mathcal{B}$ given by
$\mathcal{B}:=(A-C)\cos(\theta_0)+C\frac{\epsilon}{R}.$
 We will further manipulate
this equation to retrieve the stability results, compare also with~\cite{watanabe}. We will retrieve six groups depending on how the inertia ratio $A/C$ relates to the eccentricity $\epsilon/R$. Since in the literature results have been expressed in terms of the spin of an initial condition at a vertical state, we introduce $n_0:=\frac{J}{C(R-\epsilon)}$, which is the value of the spin $n$ at $\theta=0$ for a given Jellet $J$. Similarly, $n_{\pi}:=-\frac{J}{C(R+\epsilon)}$ is the spin of the solution with Jellet J at $\theta= \pi$. Note that for a fixed $J$ these spins are related by $n_0=-n_{\pi}\frac{R+\epsilon}{R-\epsilon}$.

\paragraph{Vertical spinning state: $\theta=0$.}

For the vertical
spinning state $\theta=0$, the relation
\eqref{stabilityba} yields
\begin{equation}\label{cond_n1}
n^2_0\left[\frac{A}{C}-(1-\frac{\epsilon}{R})\right]<
{\frac{mg\epsilon}{C}}\left(1-\frac{\epsilon}{R}\right)^2.
\end{equation}
It follows that in Group I, the vertical state
$\theta=0$ is always stable, while for Group II and Group III stability requires that
\begin{equation}\label{n1}
|n_0|<\ts n_1:=\sqrt{\frac{mg\epsilon}{C\left[\frac{A}{C}-(1-\frac{\epsilon}{R})\right]}}\left(1-\frac{\epsilon}{R}\right).
\end{equation}

\paragraph{Vertical spinning state: $\theta=\pi$.}
For the vertical
spinning state $\theta=\pi$, the relation
\eqref{stabilityba} yields
\begin{equation}\label{stabpi_muniet0}
n^2_\pi\left[(1+\frac{\epsilon}{R})-\frac{A}{C}\right]>{\frac{mg\epsilon}{C}}\left(1+\frac{\epsilon}{R}\right)^2.
\end{equation}
This condition is never satisfied for Group III, so $\theta=\pi$ is unstable; in the case of Group I and II, when $\frac{A}{C}< (1+\frac{\epsilon}{R})$, stability requires
\begin{equation}\label{n2}
|n_\pi|>\ts n_2=:\sqrt{{\frac{mg\epsilon}{C\left[(1+\frac{\epsilon}{R})-\frac{A}{C}\right]}}}\left(1+\frac{\epsilon}{R}\right)
.
\end{equation}

Note that for tippe tops of Group I and II $n_2^2\geq n_*^2$, with $n_*:=2 \frac{\sqrt{Amg\epsilon}}{C}$. The equality holds when
$\frac{A}{C}=\frac{1}{2}\left(1+\frac{\epsilon}{R}\right).$

\paragraph{Intermediate states.}
Recall that intermediate states are determined by
(\ref{define_intstates}): $f(J^2,\cos\theta_0)=0$. Using this
condition, the requirement~\eqref{stabilityba} becomes
$0<g(\cos\theta_0)$
where we set
\begin{equation}\label{eq:defg}
g(\cos\theta_0):=\ts \left(\frac{A}{C}-1\right)+\frac{4\left[(\frac{A}{C}-1)\cos\theta+\frac{\epsilon}{R}\right]^2}
{\frac{A}{C}\sin^2\theta+(\cos\theta-\frac{\epsilon}{R})^2}.
\end{equation}

We now prove that $g(\cos\theta_0)$ is strictly increasing and
changes sign at a bifurcation point for intermediate states. As we already mentioned, a bifurcation point along the curve in the
$(J^2,\cos\theta)$-plane of intermediate states is determined by the conditions
\begin{eqnarray*}
  f(J^2,\cos\theta)=&& 0,\\
  \partial_2 f(J^2,\cos\theta)=&& \ts \frac{J^2}{mg\e R^2C}(\frac{A}{C}-1) + \\
&& \ts
4((\frac{A}{C}-1)\cos\theta+\frac{\e}{R})(\frac{A}{C}(1-\cos^2\theta)
+(\cos\theta-\frac{\e}{R})^2)= 0.
\end{eqnarray*}
An elementary substitution of the first equation into the second
shows that the function $\partial_2 f(J^2,\cos\theta)$ along the
intermediate states can also be written as \[\partial_2
f(J^2,\cos\theta)=g(\cos\theta)\frac{J^2}{mg\e R^2 C}.\] Hence, bifurcation
points are given by those $\theta$ such that
$g(\cos\theta_0)=0.$
Solving for $\cos\theta$ gives the two solutions:
\begin{equation}\label{eqbif}
\frac{\e/R}{1-A/C} \pm \frac{1}{1-A/C} \frac{\sqrt3}{3}\sqrt{A/C}
\sqrt{1-A/C-(\e/R)^2}.
\end{equation}
Note that $1-A/C-(\e/R)^2>0$ only for tippe tops in Groups I and II.
Moreover, the solution $\cos\theta_b = \frac{\e/R}{1-A/C}+\ldots$
leads to a contradiction since, for tippe tops of Group I, it is
incident with the interval $]0,\theta_c[$ where no intermediate
states exist and for tippe tops in Group II satisfying
$1-A/C-(\e/R)^2>0$, the number $\frac{\e/R}{1-A/C}$ is greater than
$1$, implying that the $+$ solution of \eqref{eqbif} can not equal a
cosine. We denote the $-$ solution by $x_b$, i.e.
$$x_b : =\ts
\frac{\e/R}{1-A/C} -\frac{1}{1-A/C} \frac{\sqrt3}{3}\sqrt{A/C}
\sqrt{1-A/C-(\e/R)^2}.$$ We conclude that a bifurcation point for
intermediate states exists if $1-A/C-(\e/R)^2>0$ and
$\left|x_b\right|<1$.

Before studying in further detail these two conditions, we first
show that the function $g(\cos\theta)$ is strictly increasing for
increasing $\theta$. This result implies that, if a bifurcation
exists (i.e. a point $\theta$ with $g(\cos\theta)=0$) then stability
will change. On the other hand, if no bifurcation occurs, the entire
branch of intermediate states is either stable or unstable.

Let us assume that $x=\cos\theta$. If we consider $\partial_2
f(J^2,x)$ as a function on the submanifold $f(J^2,x)=0$ and if we
compute its derivative w.r.t $x$, i.e. $\partial_{2,2}f
-(\partial_{1,2} f)(\partial_2 f/\partial_1 f)$, then after some
tedious computations we may conclude that the sign of this
derivative is opposite to the sign of
\[\ts
8((A/C-1)x+\e/R)^2+(A/C-1)^2\frac{(A/C(1-x^2)+(x-\e/R)^2)^2}{((A/C-1)x+\e/R)^2}>0
\]
Hence, if there is a bifurcation at a certain $x=\cos(\theta_b)$ in
the set of intermediate states, then we know that the intermediate
states for which $\theta>\theta_b$ are stable, while the other
branch is unstable.

It now remains to study the conditions for the bifurcation point to
exist. The first condition says that $1-A/C-(\e/R)^2>0$, implying
that we only have to consider Groups I and II. We start with Group
I.

\paragraph{Group I} From $(A/C-1) < -\e/R$, it follows that $1-A/C-(\e/R)^2>0$ and $x_b< 1$. We have to distinguish between two subgroups: $x_b < -1$ (Group Ia) and $x_b>-1$ (Group Ib). From the previous paragraph it should be clear that if $x_b<-1$ then the
value of the function $g(x)$ on the intermediate states $-1<x<1$ is
negative, i.e. the entire branch is unstable.

\paragraph{Group II} We define three subgroups:
Group IIa is the group for which $1-A/C-(\e/R)^2>0$ and $|x_b|<1$, Group IIc is defined by $1-A/C-(\e/R)^2>0$ and $x_b< -1$ and thirdly Group IIb as the group containing the remaining tippe tops in II. Again from the previous, we
immediately conclude that the entire branch of intermediate states
is unstable in Group IIc. Group IIb can alternatively be
defined as the group containing all tippe tops for which the branch
of intermediate state is entirely stable. To show this, we first remark that the remaining tippe tops in II are characterized by $1-A/C-(\e/R)^2>0$ and $x_b>1$ or $1-A/C-(\e/R)^2\le 0$. If $1-A/C-(\e/R)^2>0$ and $x_b> 1$ then the intermediate branch is entirely stable. If $1-A/C-(\e/R)^2=0$ then $x_b>1$ and the intermediate branch is stable. The remaining tippe tops we have to consider are characterized by the condition $1-A/C-(\e/R)^2< 0$.  To show that $g(x)>0$, we consider two subcases: (i) if $1-A/C<0$ then from~\eqref{eq:defg} it is clear that $g(x) >0$, (ii) if $(\e/R)^2>1-A/C>0$ then it suffices to compute $g(0)$ ($g(x)$ will not change sign since there is no bifurcation point $x_b$): from $1-A/C-(\e/R)^2<0$,
$A/C<1$ and $(\e/R)^2<1$ we find
\[\ts
g(0)= (A/C-1) + 4\frac{(\epsilon/R)^2} {A/C+(\epsilon/R)^2}>
(\e/R)^2
>0.
\] The above argument also proves
stability for intermediate states in Group III.

\medskip
Note that tippe tops in Group II are real `tippe tops' since they
admit tipping from a position near $\theta=0$ to the inverted
state near $\theta=\pi$. Tipping never occur for tops of Group III, though they may rise up to a (stable) intermediate state. Tippe tops of group I never flip over since the position $\theta=0$ is always stable.

\subsection{Tippe Top Classification}\label{classifi}
The following schematic classification summarizes the previous analysis, and presented in Fig.~\ref{fig:bifdia}. Compare also with Fig.~3 in \cite{RW}.

\medskip \noindent
\textbf{Group I:} $A/C-1<-\e/R$ \\
- The non-inverted vertical position
$\theta=0$ is stable for any value of $J$.
\\
- The inverted
vertical position $\theta=\pi$ is stable for $|n_\pi|>n_2$, unstable otherwise, with $n_2$ given by (\ref{n2}).\\
- Intermediate states do not exist for all values of $\theta$, but only for
$\theta>\theta_c=\arccos((\frac{\epsilon}{R})/(1-\frac{A}{C}))$. \\
{\bf Group Ia}: $\ts\frac{\e/R}{1-A/C} -\frac{1}{1-A/C} \frac{\sqrt3}{3}\sqrt{A/C}
\sqrt{1-A/C-(\e/R)^2} <-1 \ .$ \\
\indent The entire branch of intermediate states
is unstable.
\\ {\bf Group Ib}: if $\ts -1<\frac{\e/R}{1-A/C} - \frac{1}{1-A/C} \frac{\sqrt3}{3}\sqrt{A/C}
\sqrt{1-A/C-(\e/R)^2}:=\cos\theta_b\ .$ \\
\indent There is a bifurcation: intermediate state are stable if  $\theta>\theta_b$ and unstable if $\theta<\theta_b$.

\medskip \noindent
\textbf{Group II:} $-\e/R<(A/C-1)<\e/R$.\\
- If $|n_0| > n_1$ the equilibria $\theta=0$ become unstable,
with $n_1$ as in (\ref{n1}).\\
- If $|n_\pi| > n_2$ the equilibria $\theta=\pi$ become
stable.\\
- There are intermediate states for any $\theta$.
We distinguish the following three subgroups. \\
{\bf Group IIa}:  $(A/C-1)<-(\e/R)^2$ and
\[\ts\left|\frac{\e/R}{1-A/C} -\frac{1}{1-A/C} \frac{\sqrt3}{3}\sqrt{A/C}
\sqrt{1-A/C-(\e/R)^2}\right|<1\ .\]
\indent There is a bifurcation of intermediate
states. \\
{\bf Group IIb}: either $(A/C-1)\ge-(\e/R)^2$ or
\[\ts \frac{\e/R}{1-A/C} -\frac{1}{1-A/C} \frac{\sqrt3}{3}\sqrt{A/C}
\sqrt{1-A/C-(\e/R)^2} >1 \ .\]
\indent The branch of intermediate states is
entirely stable.\\
{\bf Group IIc}: $(A/C-1)<-(\e/R)^2$ and
\[\ts \frac{\e/R}{1-A/C} -\frac{1}{1-A/C} \frac{\sqrt3}{3}\sqrt{A/C}
\sqrt{1-A/C-(\e/R)^2}<-1 \ .\]
\indent The branch of intermediate states is entirely unstable.

\medskip \noindent
\textbf{Group III:} $(A/C-1)>\e/R$.\\
- The equilibria with $\theta=0$ become unstable for $|n_0|>n_1$.\\
- The equilibria with $\theta=\pi$ are always unstable.\\
- For these tippe tops intermediate states do not exist for
$\theta\in ]\theta_c,\pi[$. Bifurcations in intermediate states do
not occur and the intermediate states are stable.

\medskip
\noindent
Tippe tops in Group II exhibit `tipping' behavior: if the initial spin satisfies $|n_0| > \max(n_1,n_2\frac{R+\e}{R -\e})$\footnote{About the relation between $n_1$ and
$n_2$: in Group IIb $n_2^2/n_1^2>1$ holds and in Group IIc $n_2^2/n_1^2<1$.}, then tipping is possible from $\theta=0$ to
$\theta=\pi$. 
%

\section{Remarks and Conclusions}
We would like to remark that the results on the linear stability for the tippe top presented here are equivalent to the stability properties for the model of the tippe top without the assumption that the translational friction terms can be neglected. In fact, the function $g(\cos\theta)$ which fully determines the stability of the intermediate states can be retrieved in the expressions for the eigenvalues of the linearized equations of motion about the relative equilibria of the full system. It is also remarkable that the presented stability analysis does {\em not} depend on the friction coefficient $\mu$, which suggests that the above stability analysis is valid for a rather large class of possible dissipative friction forces at the point of contact as was pointed out in~\cite{RW}. Up until now, a {\em linear stability} analysis for the intermediate states was not available upon our knowledge in the literature.

Most recently Ueda et al.~\cite{watanabe} analyzed the motion of the
tippe top under the gyroscopic balance condition ({\sc{gbc}}) and
approached the stability problem by perturbing the system around a
steady state and obtained under linear approximation a first order
{\sc{ode}} for the perturbation of the variable $\theta$. They
derive stability criteria in terms of the initial spin $n$ given at
the non-inverted position $\theta=0$. Possible intermediate states
for tippe tops of Group I were not considered. Our approach is based on a complete different technique, namely Routhian reduction and we obtain a more refined and exhaustive classification of tippe tops in six groups (instead of three). Moreover,
our analysis does not exclude the possibility of launching the top
with its stem down (i.e. $\theta$ near $\pi$).
\medskip

\section*{Acknowledgments}
The results presented here were obtained while the first author had financial support by  the European Community's 6th Framework Programme, Marie Curie Intraeuropean Fellowship EC contract Ref.\ MEIF-CT-2005-515291, award Nr.\ MATH P00286 and the second author was Postdoctoral Fellow of the Research Foundation -- Flanders (FWO) at the Department of Mathematical Physics and Astronomy, Ghent University, Belgium.
\\ The authors wish to thank Dr.\ B.\ Malengier,~Prof.\ F.\ Cantrijn, and Prof.\ J.\ Lamb for stimulating discussions and the anonymous referee for pointing out reference~\cite{RW}.


\begin{thebibliography}{99}
\bibitem{CBJB} M.C.~Ciocci, J.S.W.~Lamb, B.~Langerock and B.~Malengier. Dynamics of the spherical tippe top with small friction. {\em preprint}
\bibitem{Cohen} C.M.~Cohen. The tippe top revisited. \emph{Am. J. Phys.} 45, (1977) 12--17.
\bibitem{cush} R. Cushman. The Routh's sphere. \emph{ Reports on Mathematical Physics} 42, (1998) 47--70.
\bibitem{Duitsers} C.\ Friedl. Der Stehaufkreisel, Zulassungsarbeit zum 1.~Staatsexamen, Universit\"at Augsburg,\\
http://www.physik.uni-augsburg.de/$\sim$wobsta/tippetop/index.shtml.en

\bibitem{Gray} C.G. Gray and B.G. Nickel. Constants of the motion for nonslipping tippe tops and other tops with round pegs. \emph{Am. J. Phys.} 68 (9), (2000) 821--828.

\bibitem{Jellet} J.H. Jellet. \emph{A treatise on the theory of friction.} MacMillan, London, (1872).

\bibitem{scheck} S.\ Ebenfeld and F.\ Scheck. A new analysis of the tippe top: Asymptotic States and Lyapunov Stability. \emph{Annals of Physics} 243, (1995) 195--217.

\bibitem{Levi} T.R.~Kane and D.~Levinson. A realistic solution of the symmetric top problem. \emph{J.~Appl.~Mech.} 45, (1978) 903--909.

\bibitem{MofShi} H.K.\ Moffatt, Y.\ Shimomura and M.\ Branicki. Dynamics of axisymmetric body spinning on a horizontal surface. I. Stability and the gyroscopic approximation. \emph{Proc.\ R.\ Soc.\ Lond.} A  460, (2004) 3643--3672

\bibitem{nagrath}    I.J. Nagrath and M.~Gopal.
\newblock {\em Control Systems Engineering}.
\newblock New Age International Publishers, 4th edition, (2006).

\bibitem{nature} H.K.\ Moffatt and Y.\ Shimomura. Spinning eggs - a paradox resolved. \emph{Nature} 416, (2006) 385--386.

\bibitem{book} V.D.~Barger and M.G.~Olsson. \emph{Classical machanics. A modern perspective}. McGraw-Hill Book Company, (1973).
\bibitem{mamaev} A.V.\ Borisov and I.S.\ Mamaev. Rolling of a rigid body on plane and sphere. Hierarchy of dynamics. \emph{Regular and Chaotic Dyn.} 7 (2), (2002) 177--200.
\bibitem{Or} Or. The Dynamics of a tippe top. \emph{SIAM J.\ on A.\ Math.} 54 (3), (1994) 597--609.

\bibitem{RW} S.\ Rauch-Wojciechowski, M.\ Sk\"oldstam and T.\ Glad. Mathematical analysis of the Tippe Top. \emph{Regul. Chaotic Dyn.} 10, (2005) 333--362.

\bibitem{Romero} N.M. Bou-Rabee, J.E. Marsden and L.A. Romero. Tippe top inversion as a dissipation induced instability.
\emph{SIAM J.\ A.\ Dyn.\ Sys.} 3, (2004) 352--377.

\bibitem{routh} E.J.~Routh. \emph{Dynamics of a system of rigid bodies}. MacMillan, NY, (1905).

\bibitem{T^2} T.~Tokieda. Private Communications. The Hycaro Tipe Top was presented for the first time during the meeting \emph{Geometric Mechanics and its Applications (MASIE)}, 12--16/07/2004, EPF Lausanne, CH.
\bibitem{Tor} S.\ Torkel Glad, D.\ Petersson and S.\ Rauch-Wojciechowski. Phase Space of Rolling Solutions of the Tippe Top. \emph{ SIGMA}  3, (2007) \emph{Contribution to the Vadim Kuznetsov Memorial Issue}.

\bibitem{watanabe} T.\ Ueda, K.\ Sasaki and S.\ Watanabe. Motion of the Tippe Top: Gyroscopic balance condition and Stability. \emph{SIADS} 4 (4), (2005)  1159--1194, Society for Industrial and Applied Mathematics.

\bibitem{marsdenrouth}
J.E. Marsden, T.S. Ratiu and J.~Scheurle.
\newblock Reduction theory and the Lagrange-Routh equations.
\newblock {\em Journal of Math.\ Phys.} 41 (6), (2000) 3379--3429.

\bibitem{pars}
L.A. Pars.
\newblock {\em A Treatise on Analytical Dynamics}.
\newblock Heinemann Educational Books, (1965).



\end{thebibliography}

\appendix
\section{Routhian reduction for dissipative systems.}\label{app:routh}
Assume that a Lagrangian $L$ is given, defined on the tangent space of a manifold $M$ on which a local coordinate system $(q^1,\ldots,q^n)$ is chosen. The system is not conservative, in the sense that there is a one form on $M$, denoted by $Q=Q_idq^i$ representing the generalized force moments of the non-conservative forces. The equations of motion read:
\begin{equation}\label{thesystem2}
\frac{d}{dt}\left(\frac{\partial L}{\partial \dot q^i}
 \right)-\frac{\partial L}{\partial q^i}=Q_i \qquad \mbox{ for } i=1,\ldots,n\ .
\end{equation}
\begin{theorem}[cf.~\cite{pars}]
  Assume that there exists a coordinate, say $q^1$ such that (i) $\frac{\partial L}{\partial q^1}=0$, (ii) $Q_i(q^2,\ldots,q^n)$ is independent of $q^1$ for all $i=1,\ldots,n$, (iii) $Q_1(q^2,\ldots,q^n)=0$;
  then $\beta=\partial L/\partial \dot q^1$ is a constant of the motion. Assume that the latter equation is invertible and allows us to write $\dot q^1 =f(q^2,\ldots, q^n,\dot q^2,\ldots, \dot q^n,\beta)$. The system~\eqref{thesystem2} is equivalent to the system
  \begin{equation}\label{thesystem3}
\frac{d}{dt}\left(\frac{\partial R}{\partial \dot q^i}
 \right)-\frac{\partial R}{\partial q^i}=Q_i \qquad \mbox{ for } i=2,\ldots,n\ ,
\end{equation}
where the new Lagrangian $R$ is defined by $R= L-\dot q^1 \beta$ and such that all occurrences of $\dot q^1$ in $R$ and $Q_i$ are replaced by $f$.
\end{theorem}
By equivalent systems we mean that any solution to~\eqref{thesystem2} with fixed momentum $\partial L/\partial \dot q^1=\beta$ is a solution to~\eqref{thesystem3} and vice versa.
\section{Roots of the characteristic polynomial}\label{app:roots}
We show that the roots of the polynomial
\begin{eqnarray*} &&\lambda^3 +\mu mg\left(\frac{(R-\e\cos\theta_0)^2}{T_{\theta\theta}(\theta_0)}+ \frac{\sin^2\theta_0}{T_{\ovl\varphi\ovl\varphi}(\theta_0)}\right)\lambda^2 \\  && + \left(\frac{\big(\mu mg(R-\e\cos\theta_0)\sin\theta_0\big)^2 + \left(\fpdt{T_1}{\dot{\ovl\varphi}}{\theta}(\theta_0)\right)^2}{T_{\theta\theta}(\theta_0)T_{\ovl\varphi\ovl\varphi}(\theta_0)}  + \frac{\fpd{^2W}{\theta^2}(\theta_0)}{T_{\theta\theta}(\theta_0)}\right)\lambda
\\ && +\frac{\big(\mu mg \sin^2\theta_0\big) \fpd{^2W}{\theta^2}(\theta_0)}{T_{\theta\theta}(\theta_0)T_{\ovl\varphi\ovl\varphi}(\theta_0)}  \end{eqnarray*} all have negative real parts if and only if $\fpd{^2W}{\theta^2}(\theta_0)>0$. For our convenience we use the following shorthand notations
\begin{align*}
  &\alpha := \frac{\mu mg(R-\e\cos\theta_0)^2}{T_{\theta\theta}(\theta_0)} >0,\quad \beta:= \frac{\fpdt{T_1}{\dot{\ovl\varphi}}{\theta}(\theta_0)}{\sqrt{T_{\theta\theta}(\theta_0)T_{\ovl\varphi\ovl\varphi}(\theta_0)}},\\ & \gamma:=\frac{\mu mg\sin^2\theta_0}{T_{\ovl\varphi\ovl\varphi}(\theta_0)} > 0\mbox{ and } \delta:=\frac{\fpd{^2W}{\theta^2}(\theta_0)}{T_{\theta\theta}(\theta_0)},
\end{align*}
and the polynomial under consideration becomes \[\lambda^3 + (\alpha+\gamma)\lambda^2+ (\alpha\gamma+\beta^2+\delta)\lambda+\delta \gamma.\] Note that if $\sin\theta_0=0$ the system~\eqref{reducedeqs2} is singular and a coordinate change is necessary. The characteristic polynomial written above is not singular. Indeed, we have that $T_{\ovl\varphi\ovl\varphi}(\theta_0)$ is proportional to $\sin^2\theta_0$ or that $\gamma$ is well-defined in the case $\sin\theta_0=0$. Similarly, we have that $\fpdt{T_1}{\dot{\ovl\varphi}}{\theta}(\theta_0)$ is proportional to $\sin\theta_0$, and $\beta$ is well-defined. If we now write the necessary and sufficient conditions for this polynomial to be Hurwitz~\cite{nagrath}, we obtain the conditions
\begin{align*}
   (\alpha+\gamma)>0, (\alpha+\gamma)(\alpha\gamma+\beta^2+\delta)-\delta\gamma>0 \mbox{ and } \delta \gamma >0.
\end{align*}
The first condition is trivially satisfied. We now show that the third condition $\delta> 0$ implies the second condition. We can rewrite the second condition as
\[
(1+\alpha/\gamma)(\alpha\gamma+\beta^2+\delta) > \delta,
\]
and this is valid if $\delta>0$ or if $\fpd{^2W}{\theta^2}(\theta_0)>0$. 
\end{document}